\documentclass[notitlepage]{amsart}
\def\normalshape{\upshape}

\usepackage{amsfonts}
\usepackage[latin1]{inputenc}
\usepackage[all]{xypic}
\newtheorem{pro}{Proposition}[section]
\newtheorem{teo}[pro]{Theorem}
\newtheorem{defi}[pro]{Definition}
\newtheorem{lem}[pro]{Lemma}

\newtheorem{rk}[pro]{Remark}
\newtheorem{ex}[pro]{Example}

\newcommand{\pd}{{\mathrm{pd}}}
\newcommand{\mini}{{\mathrm{min}}}
\newcommand{\maxi}{{\mathrm{max}}}
\newcommand{\supp}{{\mathrm{Supp}_{\theta}}}

\newcommand{\pfd}{{\mathrm{pfd}}}

\newcommand{\F}{{\mathcal{F}}}
\newcommand{\I}{{\mathcal{I}}}
\newcommand{\Q}{{\mathcal{P}}}
\newcommand{\add}{{\mathrm{add}}}

\newcommand{\Ext}{\mbox{\normalshape Ext}}

\newenvironment{dem}{\noindent\bf Proof. \rm }{$\ \Box$}
\usepackage{latexsym,amssymb,amscd} 
\usepackage{amsmath}
\usepackage[all]{xy}
\begin{document}
\title{An approach to the finitistic dimension conjecture}
\author{ Fran\c cois Huard,\\ Marcelo Lanzilotta,\\ Octavio Mendoza}
\date{}
\begin{abstract} 
Let $R$ be a finite dimensional $k$-algebra over an algebraically closed
field $k$ and $\mathrm{mod}\,R$ be the category of all finitely
generated left $R$-modules. For a given full subcategory $\mathcal{X}$ of $\mathrm{mod}\,R,$ we denote by $\pfd\,\mathcal{X}$ the projective finitistic dimension of $\mathcal{X}.$ That is, $\pfd\,\mathcal{X}:=\mathrm{sup}\,\{\pd\,X\; :\; X\in\mathcal{X} \text{ and } \pd\,X<\infty\}.$
\

It was conjectured by H. Bass in the 60's that the
projective finitistic dimension $\pfd\,(R):=\pfd\,(\mathrm{mod}\,R)$ has
to be finite. Since then, much work has been done toward the proof
of this conjecture. Recently, K. Igusa and J. Todorov defined in
\cite{IT} a function $\Psi:\mathrm{mod}\,R\to \Bbb{N},$ which turned out to be useful to prove that $\pfd\,(R)$ is finite for some classes of algebras. In order to have a different approach to the finitistic dimension conjecture, we propose to consider a class of full subcategories of $\mathrm{mod}\,R$ instead of a class of algebras, namely to take the class of categories $\F(\theta)$ of $\theta$-filtered $R$-modules for all stratifying systems $(\theta,\leq)$ in $\mathrm{mod}\,R.$
\end{abstract}  

\maketitle
\section{Preliminaries}
\subsection{Basic notations}
Throughout this paper, $R$ will denote a finite dimensional $k$-algebra
over a fixed algebraically closed
field $k,$ and $\mathrm{mod}\,R$ will be the category of all finitely
generated left $R$-modules. Only finitely generated left $R$-modules
will be considered.
\

\

Given a class $\mathcal{C}$ of $R$-modules, we consider the following:
\begin{itemize}
\item[(a)] the full subcategory $\F(\mathcal{C})$ of  $\mathrm{mod}\,R$
whose objects are the zero $R$-module and the $\mathcal{C}$-filtered
$R$-modules, that is, $0\neq M \in\F(\mathcal{C})$ if there is a finite
chain $0=M_0\subseteq M_1\subseteq \cdots\subseteq M_m=M$ of
submodules of $M$ such that each quotient $M_i/M_{i-1}$ is
isomorphic to some object in $\mathcal{C};$
\item[(b)] the projective dimension
 $\pd\,\mathcal{C}:=\mathrm{sup}\,\{\pd\, C\, :\, C\in\mathcal{C}\}\in\Bbb{N}\cup\{\infty \}$ of the class $\mathcal{C};$
\item[(c)] the projective finitistic dimension
  $\pfd\,\mathcal{C}:=\mathrm{sup}\,\{\pd\, C\, :\, C\in\mathcal{C}\text{ and }
  \pd\, C<\infty\}$ of the class  $\mathcal{C};$  
\item[(d)] the class $\mathcal{P}(\mathcal{C}):=\{X\,
  :\,\Ext_R^1(X,-)|_{\F(\mathcal{C})}=0\}$ of $\mathcal{C}$-projective  $R$-modules; and 
\item[(e)] the class $\mathcal{I}(\mathcal{C}):=\{X\,
  :\,\Ext_R^1(-,X)|_{\F(\mathcal{C})}=0\}$ of $\mathcal{C}$-injective $R$-modules.
\end{itemize}
\

The global dimension $\mathrm{gldim}\,(R):=\pd\,(\mathrm{mod}\,R)$ and
the projective finitistic dimension $\pfd\,(R):=\pfd\,(\mathrm{mod}\,R)$
are important homological invariants of $R.$ The projective
finitistic dimension was introduced in the 60's  to study commutative
noetherian rings; however, it became a fundamental tool for the study of
non-commutative artinian rings. The finitistic dimension conjecture
states that $\pfd\,(A)$ is finite for any left artinian ring $A.$ This conjecture is also
closely related with other famous homological conjectures, see for
example in \cite{Xi} and \cite{ZH}.

\subsection{Stratifying systems}
For any positive integer $t\in\Bbb{Z},$ we set by definition $$[1,t]:=\{1,2,\cdots,t\}.$$
Moreover, the natural total order $\leq$ on $[1,t]$ will be considered throughout the paper. 

\begin{defi}\cite{MMS1}\label{ss}
A stratifying system  $(\theta,\leq),$ in $\mathrm{mod}\,R,$ of size $t$ consists of
a set $\theta=\{\theta(i)\, :\, i\in[1,t]\}$ of indecomposable $R$-modules  satisfying the
following homological conditions:
\begin{enumerate}
\item[(a)] $\mathrm{Hom}_R(\theta(j),\theta(i))=0$ for $j> i,$
\item[(b)] $\mathrm{Ext}^1_R(\theta(j),\theta(i))=0$ for $j\geq i.$
\end{enumerate}
\end{defi}

\subsection{Canonical stratifying systems} Once we have recalled the definition of stratifying system, it is important to know whether such a system exists for a given algebra $R.$ Consider the set $[1,n]$ which is in bijective correspondence with the iso-classes of simple $R$-modules. For each $i\in[1,n],$ we denote by $S(i)$ the simple $R$-module corresponding to $i,$ and by $P(i)$ the projective cover of $S(i).$ Following V. Dlab and C. M. Ringel in \cite{DR}, we recall that the standard module ${}_R\Delta(i)$ is the maximal quotient of $P(i)$ with composition factors amongst $S(j)$ with $j\leq i.$ So, by Lemma 1.2 and Lemma 1.3 in \cite{DR}, we get that the set of standard modules ${}_R\Delta:=\{{}_R\Delta(i)\, :\,i\in[1,n]\}$ satisfies Definition \ref{ss}. We will refer to this set as the cannonical stratifying system of $R$. Moreover, in case ${}_RR\in\F({}_R\Delta)$ following I. Agoston,V. Dlab and E. Lukacs in \cite{ADL}, we say that $R$ is a standardly stratified algebra (or a ss-algebra for short).

\subsection{The finitistic dimension conjecture for ss-algebras} It was shown in \cite{AHLU} that the finitistic dimension conjecture holds for ss-algebras. In this paper, it is shown that $\pfd\,(R)\leq 2n-2$ where $n$ is the number of iso-classes of simple $R$-modules. Of course, not only the number $n$ but also the category $\F({}_R\Delta)$ is closely related to $\pfd\,(R)$ as can be seen in the following result (see the proof in the Corollary 6.17(j) in \cite{MS}).
\

\begin{teo}\cite{MS} \label{MS1} Let $R$ be a ss-algebra and $n$ be the number of iso-classes of simple $R$-modules. Then
$$\pfd\,(R)\leq \pd\,\F({}_R\Delta)+\mathrm{resdim}_{\F({}_R\Delta)}\,(\F({}_R\Delta)^\wedge)\leq 2n-2.$$
\end{teo}

\section{The stratifying system approach}

In Theorem \ref{MS1}, we saw that if $R$ is a ss-algebra, then the category $\F({}_R\Delta)$ is closely related to $\pfd\,(R).$ Assume now that $R$ is not a ss-algebra. We know from 1.3 that $R$ still admits at least one stratifying system. Thus it makes sense to study $\pfd\,(R)$ from the point of view of stratifying systems. The following theorem was shown in \cite{MS}: 
\

\begin{teo}\cite{MS} \label{MS2} Let $R$ be an algebra, $n$ be the number of iso-classes of simple $R$-modules and $(\theta,\leq)$ be a stratifying system of size $t.$ If $\mathcal{I}(\theta)$ is coresolving, then
$$\pd\,\F(\theta)\leq t\leq n\; \text{ and }\quad \pfd\,(R)\leq \pd\,\F(\theta)+\mathrm{resdim}_{\mathcal{Y}}\,(\mathcal{Y}^\wedge),$$ where $\mathcal{Y}:=\{X\, :\,\Ext_R^1(X,-)|_{\mathcal{I}(\theta)}=0\}\supseteq \F(\theta).$
\end{teo}

\begin{rk} (1) The condition of $\mathcal{I}(\theta)$ being coresolving given in the Theorem \ref{MS2} is strong and has as a consequence that $\pd\,\F(\theta)$ is finite. However, we do not know whether the resolution dimension $\mathrm{resdim}_{\mathcal{Y}}\,(\mathcal{Y}^\wedge)$ is finite.
\

(2) In general, we can have that $\pd\,\F(\theta)=\infty$ as can be seen in the following example.
\end{rk}

\begin{ex} Let $R$ be the quotient path algebra given by the quiver
$$
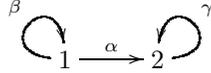
\xymatrix{\ar@(l,u)[]^{\beta} 1 \ar[r]^{\alpha} & 2\ar@(r,u)[]_{\gamma}}$$
and the relations $\beta^2=\alpha\beta=\gamma\alpha=\gamma^2=0.$ We have that
${}_R\Delta(1)=P(1)/S(2)$ and ${}_R\Delta(2)=P(2).$ Then
$({}_R\Delta,\leq)$ is a stratifying system of size $2$ and $\pd\,{}_R\Delta(1)=\infty.$
\end{ex}
\

 The discussion above leads us to the following questions. Question 1: suppose that $\pd\,\F(\theta)=\infty;$ is it possible to prove that $\pfd\,\F(\theta)<\infty$? Of course, if the finitistic dimension conjecture is true, we have that $\pfd\,\F(\theta)$ is finite.
\

Question 2: is the following number finite?
$$\mathrm{sspfd}\,(R):=\mathrm{sup}\,\{\pfd\,\F(\theta) \,:\, (\theta,\leq)\text{ is a stratifying system with } \pfd\,\F(\theta)<\infty\}.$$
Again we have that if $\pfd\,(R)$ is finite then so is $\mathrm{sspfd}\,(R)$. Hence, we propose the following questions.
\

\

{\bf The Homological Stratifying System Conjecture I:} {\it Let $R$ be any finite dimensional $k$-algebra. Then, for any stratifying system $(\theta,\leq)$ in $\mathrm{mod}\,R,$ we have that $\pfd\,\F(\theta)$ is finite}.
\

{\bf The Homological Stratifying System Conjecture II:} {\it For any finite dimensional $k$-algebra $R,$ we have that $\mathrm{sspfd}\,(R)$ is finite.}
\

\

Note that the concept of stratifying system could be an appropriate tool to
construct a counter-example for the finitistic dimension conjecture
since $\pfd\,\F(\theta)\leq \pfd\,(R).$

\begin{rk}\label{trivial} Let $R$ be a $k$-algebra and $(\theta,\leq)$ be a stratifying
system of size $t.$ If one of the following three conditions holds: (a) t=1, (b) $\pd\,\theta<\infty$ or (c) $\I(\theta)$ is coresolving, then  $\pfd\,\F(\theta)$ is finite. Indeed, if  If $t=1,$ we have that $\F(\theta)=\add\,\theta(1)$ with $\theta(1)$ indecomposable. On the other hand, we know, by the Corollary 2.4 in \cite{MMS3}, that $\pd\,\F(\theta)=\pd\,\theta.$ Finally, if $\I(\theta)$ is coresolving, we obtain from \ref{MS2} that $\pfd\,\F(\theta)$ has to be finite.
\end{rk} 

\section{Igusa-Todorov's function}
In this section, we recall the Igusa-Todorov function \cite{IT} and some properties of this function that have been successfully applied to solve the finitistic dimension conjecture for some special classes of algebras.  Moreover, we establish a new interesting inequality (see \ref{LemaH}), which is our contribution, and will be useful to give a partial solution to The Homological Stratifying System Conjecture I.
\

Let $\mathcal{K}$ denote the quotient of the free abelian group generated by all 
the symbols $[M],$ where $M\in\mathrm{mod}\,R,$ modulo the relations: (a) $[C]=[A]+[B]$ if
$C\simeq A\coprod B,$ and (b) $[P]=0$ if $P$ is projective. Then
$\mathcal{K}$ is the free $\Bbb{Z}$-module generated by the iso-classes
of indecomposable finitely generated  non-projective $R$-modules. In
\cite{IT}, K. Igusa and G. Todorov defined a function
$\Psi:\mathrm{mod}\,R\to \Bbb{N}$ as follows. Denote the first syzygy of $M$ by $\Omega\,M.$ 
\

The syzygy induces a $\Bbb{Z}$-endomorphism on $\mathcal{K}$ that will also be 
denoted by $\Omega.$ That is, we have a $\Bbb{Z}$-homomorphism  
$\Omega :\mathcal{K}\rightarrow \mathcal{K}$ where $\Omega [M]:=[\Omega\,M].$ 
For a given $R$-module $M,$ we denote by $<M>$ the $\Bbb{Z}$-submodule of $\mathcal{K}$ generated by the indecomposable direct summands of $M.$  Since $\Bbb{Z}$ is a Noetherian ring, Fitting's lemma implies that there is an integer $n$ such that  $\Omega:\Omega^m<M>\;\rightarrow\Omega^{m+1}<M>$ is an isomorphism for all $m\geq n;$ hence there exists
a smallest non-negative integer $\Phi\,(M)$ such that $\Omega:\Omega^m<M>\;\rightarrow
\Omega^{m+1}<M>$ is an isomorphism for all $m\geq \Phi\,(M).$ Define
$\mathcal{C}_M$ as the set whose elements are the direct summands of
$\Omega^{\Phi(M)}\,(M).$ Then we set:
$$\Psi\,(M):=\Phi\,(M)+\pfd\,\mathcal{C}_M.$$ The following result is due to K. Igusa and
G. Todorov.

\begin{pro}\cite{IT}\label{ProIT} The function $\Psi:\mathrm{mod}\,R\to\Bbb{N}$ satisfies the following properties.
\begin{itemize}
\item[(a)] If $\pd\, M$ is finite then $\Psi\,(M)=\pd\,M.$ On the other hand, $\Psi\,(M)=0$ if $M$ is indecomposable and $\pd\,M=\infty.$
\item[(b)] $\Psi\,(M)=\Psi\,(N)$ if $\add\,M=\add\,N.$
\item[(c)] $\Psi(M)\leq\Psi(M\coprod N).$ 
\item[(d)] $\Psi(M\coprod P)=\Psi(M)$ for any projective $R$-module $P.$
\item[(e)] If $0\to A\to B\to C\to 0$ is an exact sequence in $\mathrm{mod}\,R$ and $\pd\,C$ is finite then $\pd\,C\leq\Psi\,(A\coprod B)+1.$  
\end{itemize}
\end{pro}
\

The following useful inequality was first proved by Y. Wang.

\begin{lem}\label{W} \cite{W} If $0\to A\to B\to C\to 0$ is an exact sequence in $\mathrm{mod}\,R$ and $\pd\,B$ is finite then $\pd\,B\leq 2+\Psi\,(\Omega\,A\coprod \Omega^2\,C).$
\end{lem}
\

As can be seen above, the Igusa-Todorov function $\Psi$ is some kind of generalization of the projective dimension. So, it makes sense to consider the ``$\Psi$-dimension" of a given class $\mathcal{X}$ of objects in $\mathrm{mod}\,R.$

\begin{defi}\label{quest} Let $\mathcal{X}$ be a class of objects in $\mathrm{mod}\,R.$ The $\Psi$-dimension of $\mathcal{X}$ is $\Psi\mathrm{d}\,(\mathcal{X}):=\mathrm{sup}\,\{\Psi (X)\,:\,X\in\mathcal{X}\}.$ 
\end{defi}

\noindent {\bf{Question}}: Is $\Psi\mathrm{d}\,(\F(\theta))$ finite for any stratifying system $(\theta,\leq)$ in $\mathrm{mod}\,R$?
\

\

It would be useful to have the following generalisation of Proposition \ref{ProIT} (e): if $0\to A\to B\to C\to 0$ is an exact sequence in $\mathrm{mod}\,R$ then $\Psi\,(C)\leq \Psi\,(A\coprod B)+1.$  Such a result is not true in general but we will now show that it holds when $B$ is projective. It will be an appropriate tool to obtain a bound of $\pfd\,\mathcal{X}$ for a given class $\mathcal{X}$ of objects in $\mathrm{mod}\,R.$
\

\begin{lem}\label{paraLemaH} 
$\Phi\,(X)\leq 1+ \Phi\,(\Omega\,X)$ for any $X\in\mathrm{mod}\,R.$
\end{lem}
\begin{dem} Let $X$ be an $R$-module; we assert that $\Omega <X>$ is a
  subgroup of $<\Omega X>.$ Indeed, let $X_1,X_2,\cdots,X_t$ be all
  the indecomposable (up to isomorphism) non projective direct
  summands of $X.$ For each $i,$ we have a decomposition
  $\Omega\,X_i=\oplus_{k=1}^{n_i}X_{ik}^{\beta_{ik}}\oplus P_i$  where
  $X_{ik}$ are pairwise non-isomorphic and $P_i$ is a projective $R$-module. Let $x\in <X>,$ so we have $x=\sum_{i=1}^t\alpha_i[X_i]$ with $\alpha_i\in\Bbb{Z}$ for each $i.$  Applying $\Omega$ to $x,$ we get the following equalities
$$\Omega\,(x)=\sum_{i=1}^t\alpha_i[\Omega\,X_i]=\sum_{i=1}^t\alpha_i(\sum_{k=1}^{n_i}\beta_{ik}[X_{ik}]);$$ proving that  $\Omega <X>\subseteq <\Omega X>.$ 
\

Assume that $\Phi\,(\Omega\,X)=n.$ Hence, $\Omega:\Omega^m<\Omega\,X>\;\rightarrow \Omega^{m+1}<\Omega\,X>$ is an isomorphism for all $m\geq n.$ On the other hand, $\Omega^{m+1}<X>$ is a subgroup of $\Omega^m<\Omega\,X>$ since $\Omega <X>\subseteq <\Omega X>.$ Therefore, we get an isomorphism $\Omega:\Omega^{m+1}<X>\;\rightarrow \Omega^{m+2}<X>$ for all $m+1\geq
n+1.$ Then $\Phi\,(X)\leq n+1=\Phi\,(\Omega\,X)+1.$
\end{dem}

\begin{pro}\label{LemaH} $\Psi\,(X)\leq 1+ \Psi\,(\Omega\,X)$ for any $X\in\mathrm{mod}\,R.$
\end{pro}
\begin{dem} Let $X$ be an $R$-module. Consider the set $\mathcal{C}_X$ whose elements are the direct summands of $\Omega^{\Phi\,(X)}\,(X);$ and choose $Y$ in $\mathcal{C}_X$ such that $\pd\,Y=\pfd\,\mathcal{C}_X.$ So, by the definition of $\Psi,$ we have that $\Psi\,(X)=\pd\,Y+\Phi\,(X).$ Hence, to prove the result, it is enough to see that $\pd\,Y+\Phi\,(X)-1\leq \Psi\,(\Omega\,X).$
\

Take $n:=\Phi\,(X)-1$ and $M:=\Omega\,X.$ Then, we have that $Y$ is a direct summand of $\Omega^nM$ since $\Omega^n\,M=\Omega^{\Phi\,(X)}\,(X).$ Therefore, using that $n\leq\Phi\,(M)$(see \ref{paraLemaH}), we get from Lemma 3 d) in \cite{IT} that $\pd\,Y+n\leq \Psi\,(M),$ proving the result.
\end{dem}

\begin{lem}\label{sisieq} Let $0\to A\to B\to C\to 0$ be an exact sequence of $R$-modules. Then we have the following:
\begin{itemize}
\item[(a)] if $\pd\,C$ is finite then for any $m\geq \pd\,C$ there are projective $R$-modules $P_m$ and $P'_m$ such that $\Omega^{m}\,(A)\coprod P_m\simeq\Omega^{m}\,(B)\coprod P'_m,$ 
\item[(b)] if $\pd\,A$ is finite then for any $m\geq \pd\,A$ there are  projective $R$-modules $P_m$ and $P'_m$ such that $\Omega^{m+1}\,(B)\coprod P_m\simeq\Omega^{m+1}\,(C)\coprod P'_m,$ 
\item[(c)] if $\pd\,B$ is finite then for any $m\geq \pd\,B+1$ there are projective $R$-modules $P_m$ and $P'_m$ such that $\Omega^{m+1}\,(C)\coprod P_m\simeq\Omega^{m}\,(A)\coprod P'_m.$ 
\end{itemize}
\end{lem}
\begin{dem} The proof follows from the well known fact: for any exact sequence $0\to X\to Y\to Z\to 0$ of $R$-modules, there are two new exact sequences $0\to\Omega\,Y\to \Omega\,Z\coprod P\to X\to 0$ and  $0\to \Omega\,X\to \Omega\,Y\coprod P'\to \Omega Z\to 0,$ where $P$ and $P'$are projective $R$-modules. 
\end{dem}

\section{Main results}

In this section, we give a partial solution to The Homological Stratifying System Conjecture I, which was stated in Section 2. In order to do that, we have to introduce some definitions and also to recall several properties for stratifying systems that were proven in \cite{MMS2}.
\

 Let $(\theta,\leq)$ be a stratifying system of size $t.$ Due to K. Erdmann and C. S\'aenz in \cite{ES}, we know that the filtration multiplicities $[M:\theta(i)]$ do not depend on the filtration of $M\in\mathcal{F}(\theta).$ Following the notation in \cite{MMS2}, we recall that the $\theta$-support of $M$ is the set $\mathrm{Supp}_{\theta}(M)=\{i\in[1,t]\; :\; [M:\theta(i)]\not=0\}.$ Therefore, $\mathrm{Supp}_{\theta}(M)$ is empty if $M=0.$ The functions $\mathrm{min},\;\mathrm{max}:\F(\theta)\rightarrow [1,t]\cup\{\pm\infty\}$ are defined as follows: (a) $\mathrm{min}\,(0):=+\infty$ and $\mathrm{max}\,(0):=-\infty,$ and (b) $\mathrm{min}(M):=\mathrm{min}\,(\mathrm{Supp}_{\theta}(M),\leq)$ and $\mathrm{max}(M):=\mathrm{max}\,(\mathrm{Supp}_{\theta}(M),\leq)$ if $M\neq 0.$  For a given $0\neq M\in\F(\theta),$ we will denote by $(\theta_M,\leq)$ the new stratifying system induced by the set $\theta_M:=\{\theta(i)\,:\,i\in\supp\,(M)\}.$ Such kind of reduction will be very useful since $M\in\F(\theta_M)$ and the size of $(\theta_M,\leq)$ is smaller than $t.$

\begin{defi}\label{epss} \cite{MMS2} 
Let $\theta=\{\theta(i)\,:\,i\in[1,t]\}$ be a set of non-zero $R$-modules and $\underline{Q}=\{Q(i)\,:\,i\in[1,t]\}$ be a set of indecomposable  $R$-modules. We say that the system $(\theta,\underline{Q}, \leq )$ is an Ext-projective stratifying system ( or epss for short) of size $t$, if the following three conditions hold:
\begin{enumerate}  
\item $\mathrm{Hom}_R(\theta(j),\theta(i))=0$ for $j>i,$
\item for each $i\in[1,t]$ there is an exact sequence $0\to K(i)\rightarrow Q(i)\rightarrow\theta(i) \to 0$ such that  $K(i)\in\mathcal{F}(\{\theta(j): j>i\}),$
\item $\mathrm{Ext}_R^1(Q,\mathcal{F}(\theta))=0$, where $Q:=\coprod_{i=1}^tQ(i).$ 
\end{enumerate}
\end{defi}

\begin{rk} Due to \cite{MMS2}, we have that, for a given stratifying system $(\theta,\leq),$ there is a unique - up to isomorphism -  Ext-projective stratifying system $(\theta,\underline{Q},\leq)$ which is called the epss associated to $(\theta,\leq).$
\end{rk}
\

The following results, which were shown in \cite{MMS2}, will be very useful throughout this section.

\begin{pro}\label{suc} \cite{MMS2} Let $(\theta,\leq)$ be a stratifying system of size $t,$ $0\not=M\in\mathcal{F}(\theta),$ $i=\mathrm{min}(M)$ and $m=[M:\theta(i)].$ Then:
\begin{enumerate}
\item there exists a finite chain
$$ 0\subset N\subset M_{m-1}\subset \cdots\subset M_1\subset M_0=M,$$
of submodules of $M$ sucht that $N\in \mathcal{F}(\{\theta(j):j>i \})$ and $M_k/M_{k+1}\cong \theta(i)$ for all $k=0,1,\cdots, m-1,$ where $M_m=N$.
\item there exists an exact sequence in $\F(\theta)$ $$0\rightarrow N\rightarrow
M\rightarrow \theta(i)^{m}\rightarrow 0 \quad\text{ with }\quad\mathrm{min}(M)<\mathrm{min}(N).$$
\end{enumerate}
\end{pro}
\

Given a stratifying system $(\theta,\leq),$ the category $\F(\theta)$ of $\theta$-filtered $R$-modules has very nice properties; for example, it has Ext-projective covers. We recall that, a morphism $f:Q_0(M)\rightarrow M$ in $\F(\theta)$ is an Ext-projective cover of $M$ if and only if $f:Q_0(M)\rightarrow M$ is the right minimal $\mathcal{P}(\theta)$-approximation of $M$ and $\mathrm{Ker}\,f\in\F(\theta).$ The following result says in particular that any $M\in\F(\theta)$ has an Ext-projective cover. In fact, following Definition 4.1, we have that $Q_0(\theta(i))=Q(i)$ for any $i;$ and so, $Q=Q_0(\coprod_{i=1}^t \theta(i)).$

\begin{pro}\label{sucepss} \cite{MMS2} Let $(\theta,\leq)$ be a stratifying system of size $t$ and $(\theta,\underline{Q},\leq)$ be the epss associated to $(\theta,\leq).$ Let  $0\not=M\in
\mathcal{F}(\theta),$ $i=\mathrm{min}(M)$ and $m_i=[M:\theta(i)].$
Then, there exists an exact sequence in $\F(\theta)$ $$0\rightarrow N\rightarrow
Q_0\,(M)\stackrel{\varepsilon_M}{\rightarrow} M\rightarrow 0 $$ such that $\mathrm{min}(M)<\mathrm{min}(N),$  $Q_0\,(M)\in\add\,\coprod_{j\geq i}\,Q(j)$ and $\varepsilon_M\,:Q_0\,(M)\rightarrow M$ is the right minimal $\Q(\theta)$-approximation of $M.$ 
\end{pro}

\begin{defi}\label{inftheta} Let $(\theta,\leq)$ be a stratifying system of size $t.$  We denote by $\infty_\theta$ the set $\{i\in[1,t]\,:\, \pd\,\theta(i)=\infty\},$ and by $\mathrm{card}\,\infty_\theta$ the cardinality of $\infty_\theta.$
\end{defi}

\begin{lem}\label{LemaH2} Let $(\theta,\leq)$ be a stratifying system of size $t$ such that $\infty_\theta=\{i_0\}.$ Define $s:=\mathrm{max}\,\{\pd\,\theta(j)\;:\; j\neq i_0\}$ if $t>1,$ and $s:=0$ otherwise. Then, for any $M\in\mathcal{F}(\theta),$ we have that
\begin{itemize}
\item[(a)] $\pd\,M<\infty$ if and only if $i_0\not\in\supp\,(M),$ 
\item[(b)] if $i_0\in\supp\,(M),$ then there exist projective $R$-modules $P$ and $P'$ such that $$\Omega^{s+1}\,(M)\coprod P\simeq\Omega^{s+1}(\theta(i_0)^{[M:\theta(i_0)]})\coprod P'.$$ 
\end{itemize}
\end{lem}
\begin{dem} If $t=1$ there is nothing to prove since we know that $\mathcal{F}(\theta)=\add\,\theta(i_0)$ and $\pd\,0=0.$ So, we may assume that $t>1.$
\

(a) $(\Leftarrow)$ If $i_0\not\in\supp\,(M),$ then by considering the new stratifying system induced by the set $\theta_M:=\{\theta(i)\; :\;
i\in\mathrm{Supp}_\theta\,M\},$ we get from \ref{trivial} that $\pd\,M\leq\pd\,\F(\theta_M)=\pd\theta_M\leq s.$
\

$\qquad(\Rightarrow)$ Let $M\in\mathcal{F}(\theta)$ be of finite projective dimension. Suppose that $[M:\theta(i_0)]=m\neq 0.$ Then, by \ref{suc}, we have a chain $M_{i_0}\subseteq M_{i'_0}\subseteq M$ of submodules of $M$ such that: $M_{i'_0}/M_{i_0}\simeq \theta(i_0)^m,$  $i_0\not\in\supp\,(M_{i_0})$ and $i_0\not\in\supp\,(M/M_{i'_0}).$ Consider the following exact sequences in $\F(\theta)$
\begin{equation}\label{S2f1} 0\to M_{i_0}\to M\to M/M_{i_0}\to 0,
\end{equation}
\begin{equation}\label{S2f2} 0\to M_{i'_0}/M_{i_0}\to M/M_{i_0}\to M/M_{i'_0}\to 0.
\end{equation}
Because of the fact that $i_0\not\in\supp\,(M_{i_0}),$ we get that $\pd\,M_{i_0}$ is finite (see the proof of $(\Leftarrow)$ given above). Hence, by \eqref{S2f1}, we obtain that $\pd\,M/M_{i_0}$ is finite.\\
We assert that $M_{i'_0}\neq M.$ Otherwise, we would obtain that $M/M_{i_0}\simeq\theta(i_0)^m,$ and so $\pd\,M/M_{i_0}$ has to be infinite which is a contradiction, proving that  $M_{i'_0}\neq M.$\\
Since $M_{i'_0}\neq M$ and $i_0\not\in\supp\,(M/M_{i'_0}),$ we have that $\pd\,M/M_{i'_0}$ is finite; consequently, by \eqref{S2f2}, we conclude that $\pd\,M_{i'_0}/M_{i_0}$ is finite which is a contradiction, proving that $i_0\not\in\supp\,(M).$ 
\

(b) Assume that $i_0\in\supp\,(M).$ Let $m:=[M:\theta(i_0)];$ so, from \ref{suc}, we get the following two exact sequences
\begin{equation}\label{S2f3} 
0\to M_0\to M\to M/M_0\to 0,\quad 0\to N\to M_0\to \theta(i_0)^m\to 0
\end{equation}
such that $\mathrm{min}\,(M_0)=i_0$ and $[M/M_0:\theta(i_0)]=0=[N:\theta(i_0)].$ Hence, as we have proven in (a), we have that $\pd\,M/M_0\leq s$ and $\pd\,N\leq s.$ Therefore, by applying \ref{sisieq} to the exact sequences in \eqref{S2f3}, we conclude that $\Omega^{s+1}(M_0)\simeq\Omega^{s+1}(M)$ and the existence of projective $R$-modules $P$ and $P'$ such that $\Omega^{s+1}(M_0)\coprod P\simeq\Omega^{s+1}(\theta(i_0)^m)\coprod P';$ proving the result.
\end{dem}

\begin{teo}\label{Teofd1}  Let $(\theta,\leq)$ be a stratifying system of size $t$ such that $\infty_\theta=\{i_0\}.$ Define $s:=\mathrm{max}\,\{\pd\,\theta(j)\;:\; j\neq i_0\}$ if $t>1,$ and $s:=0$ otherwise. Then
\begin{itemize}
\item[(a)] $\pfd\,\F(\theta)\leq s,$
\item[(b)] $\Psi\mathrm{d}\,(\F(\theta))\leq 1+s+\Psi\,(\Omega^{s+1}\,\theta(i_0)).$
\end{itemize}
\end{teo}
\begin{dem} If $t=1,$ we have that $\mathcal{F}(\theta)=\add\,\theta(i_0);$ and so, by \ref{ProIT}, we conclude that $\pfd\,\F(\theta)=0=\Psi\mathrm{d}\,(\F(\theta)).$ Hence, we may assume that $t>1.$
\

(a)  Let $M\in\F(\theta)$ be of finite projective dimension. Consider the new stratifying system $\theta_M:=\{\theta(i)\; :\; i\in\mathrm{Supp}_\theta\,M\}$ induced by $M.$ It follows, from \ref{LemaH2} (a), that $i_0\not\in\mathrm{Supp}_\theta\,M.$ Hence $\pd\,M\leq \pd\,\theta_M\leq s.$
\

(b) Let $ M\in\F(\theta)$ be such that $\pd\,M=\infty.$ Then by \ref{LemaH2} (a) we have that $i_0\in\supp\,(M).$ Hence, from \ref{LemaH2} (b) and \ref{ProIT}, we get that $\Psi\,(\Omega^{s+1}M)=\Psi\,(\Omega^{s+1}\theta(i_0)).$ Therefore, using \ref{LemaH}, we conclude that $\Psi\,(M)\leq 1+s+\Psi\,(\Omega^{s+1}\theta(i_0)).$
\end{dem} 
\

Note that the previous result gives a partial answer to the Question \ref{quest}; and as a consequence, we get the first Homological Stratifying System Conjecture in this case. We will now show that the latter holds in case $(\theta, \leq)$ admits at most two indecomposable modules of infinite projective dimension. 

\begin{teo}\label{Teofd2}  Let $(\theta,\leq)$ be a stratifying system of size $t$ such that $\infty_\theta=\{i_0,i_1\},$ where $i_0<i_1.$ Define  $s:=\mathrm{max}\,\{\pd\,\theta(j)\;:\; j\neq i_0,i_1\}$ if $t>2,$ and $s:=0$ otherwise. Then, we have that 

$$\pfd\,\F(\theta)\leq s+2+\mathrm{min}\,(\alpha,\beta)$$ 
\

\noindent where $\alpha:=\Psi\,(\Omega^{s+1}\,\theta(i_1)\coprod\Omega^{s+2}\,\theta(i_0)),$ 
$\beta:=\Psi(\Omega^{s+1}\,(Q\coprod\theta(i_1)))$ and $Q$ is the Ext-projective cover of 
$\coprod_{i=1}^t\theta(i)$ in $\F(\theta).$
\end{teo}
\begin{dem} Take $M\in\F(\theta)$ with $\pd\,M<\infty.$ Let $j\in\infty_\theta,$ if $j\not\in\supp\,(M)$ then by applying \ref{Teofd1} (a) to  $(\theta_M,\leq),$ we get that $\pd\,M\leq s.$ So, we may assume that $i_0,i_1\in\supp\,(M).$ By \ref{suc}, we have a $\theta$-filtration of $M$
$$M_1\subseteq M'_{i_1}\subseteq\cdots \subseteq M_{i_1}\subseteq\cdots\subseteq M_0\subseteq
M'_{i_0}\subseteq\cdots\subseteq M_{i_0}\subseteq\cdots\subseteq M $$
such that $[M:\theta(i_0)]=[M_{i_0}/M_0:\theta(i_0)]$ and  $[M:\theta(i_1)]=[M_{i_1}/M_1:\theta(i_1)].$ Consider the exact sequence $0\to M_{i_0}\to M\to M/M_{i_0}\to 0.$ Since $i_0,i_1\not\in\supp\,(M/M_{i_0}),$ we get $\pd\,M/M_{i_0}\leq s;$  therefore $\pd\,M\leq\maxi\,(s,\pd\,M_{i_0})$ and $\pd\,M_{i_0}<\infty.$ On the other hand, from the exact sequence $0\to M_1\to M_{i_0}\to M'\to 0,$ we obtain that $\pd\,M_{i_0}\leq\maxi\,(s,\pd\,M'),$ and also that 
$\pd\,M'<\infty$ since $\pd\,M_1$ is finite ($i_0,i_1\not\in\supp\,(M_1)$). Hence $\pd\,M\leq\maxi\,(s,\pd\,M')$ with $\pd\,M'$ finite, $\mini\,(M')=i_0$ and $\maxi\,(M')=i_1.$ 
\

We assert that $\pfd\,\F(\theta)\leq s+2+\alpha.$ Indeed, let $m_1:=[M':\theta(i_1)],$ then we have an exact sequence in $\F(\theta)$ 
\begin{equation}\label{f1} 0\rightarrow\theta(i_1)^{m_1}\rightarrow
M'\rightarrow M''\rightarrow 0 
\end{equation}
such that $\supp\,(M'')\cap\infty_{\theta}=\{i_0\}.$ Since $\pd\,M'$ is finite, we get, from \ref{W}, that $\pd\,M'\leq 2+\Psi\,(\Omega\,\theta(i_1)\coprod \Omega^2\,M'').$ On the other hand, using that $\supp\,(M'')\cap\infty_{\theta}=\{i_0\},$ we get from \ref{LemaH2} that $\Omega^{s+2}\,M''\simeq \Omega^{s+2}\,(\theta(i_0)^{[M:\theta(i_0)]});$ and so, from \ref{LemaH}, we conclude that $\pd\,M\leq s+2+\alpha;$ proving our assertion.
\

Finally, we prove that  $\pfd\,\F(\theta)\leq s+2+\beta.$ In order to do that, we consider the epss $(\theta,\underline{Q},\leq)$ associated to $(\theta,\leq)$. So, as we have recalled before, we have that $Q:=\coprod_{i=1}^t Q(i)$ is the Ext-projective cover of $\coprod_{i=1}^t \theta(i)$ in $\F(\theta)$. By \ref{sucepss}, we obtain the following exact sequence in $\F(\theta)$
\begin{equation}\label{f2} 0\rightarrow N\rightarrow
Q_0\,(M')\rightarrow M'\rightarrow 0, 
\end{equation}
where $i_0=\mini(M')<\mathrm{min}(N)$ and
$Q_0\,(M')\in\add\,\coprod_{j\geq i_0}\,Q(j).$ In particular, we have
that $i_0\not\in\supp\,(N).$ If
$\supp\,(N)\cap\infty_\theta=\emptyset,$ we have that $\pd\,N\leq s.$
Hence, we can apply \ref{sisieq} (b) to \eqref{f2} obtaining that $\Psi\,(\Omega^{s+1}\,Q_0(M'))=\Psi\,(\Omega^{s+1}\,M').$ Therefore, from \ref{LemaH}, we get $\pd\,M'=\Psi\,(M')\leq s+1+\Psi\,(\Omega^{s+1}\,Q_0(M'))\leq s+1+\Psi\,(\Omega^{s+1}\,Q)\leq s+2+\beta$ since $Q_0(M')\in\add\,Q.$\\
\

Suppose that $\supp\,(N)\cap\infty_\theta=\{i_1\}.$ Then, applying \ref{LemaH2} (b) to the stratifying system $(\theta_N,\leq),$ we get that $\Omega^{s+1}\,(N)\coprod P\simeq \Omega^{s+1}\,(\theta(i_1)^{[N:\theta(i_1)]})\coprod P',$ where $P$ and $P'$ are projective $R$-modules. Hence, from \ref{ProIT} and \ref{LemaH}, we conclude that $\Psi\,(N\coprod Q_0(M))\leq s+1+\Psi\,(\Omega^{s+1}\,(Q\coprod\theta(i_1))).$ On the other hand, applying \ref{ProIT} (e) to \eqref{f2}, we have that $\pd\,M\leq 1+\Psi\,(N\coprod Q_0(M)).$ Therefore $\pd\,M\leq s+2+\Psi\,(\Omega^{s+1}\,(Q\coprod\theta(i_1));$ proving the result.
\end{dem}
\

\

Let $(\theta,\leq)$ be a stratifying system. As can be seen above, we have proven that $\pfd\,\F(\theta)$ is finite in the cases $\mathrm{card}\,\infty_{\theta}=1,2$ without extra conditions on $\F(\theta).$ However, if $\mathrm{card}\,\infty_{\theta}=3,$ we were not able to prove that $\pfd\,\F(\theta)$ is finite without extra conditions on $\F(\theta).$ In order to get a proof, we will assume one of the following ``3-properties".
\

\begin{defi}\label{3-properties} Let $(\theta,\leq)$ be a stratifying system of size $t$ with $\mathrm{card}\,\infty_{\theta}=3.$ That is, there are indices $i_0<i_1<i_2$  such that $\infty_\theta=\{i_0,i_1,i_2\}.$ Let $M$ be any element of $\F(\theta)$ and $\varepsilon_M: Q_0(M)\rightarrow M$ be the Ext-projective cover of $M.$ We say that
\

(a) $\F(\theta)$ satisfies the {\bf{ 3-finitistic property}} if $i_1,i_2\in\supp\,(\mathrm{Ker}\,\varepsilon_M)$ and $\pd\,M<\infty$ implies that $\pd\,\mathrm{Ker}\,\varepsilon_M$ is finite; and
\

(b) $\F(\theta)$ satisfies the {\bf{ 3-cardinal property}} if $\supp\,(M)\cap\infty_\theta=\{i_0,i_1\}$ implies that $\mathrm{card}\,(\supp\,(\mathrm{Ker}\,\varepsilon_M)\cap\infty_\theta)\leq 1.$ 
\end{defi}

\begin{teo}\label{Teofd3}  Let $(\theta,\leq)$ be a stratifying system of size $t$ such that $\infty_{\theta}=\{i_0,i_1,i_2\}$ and $i_0<i_1<i_2.$ Consider  $s:=\mathrm{max}\,\{\pd\,\theta(j)\;:\; j\not\in\infty_\theta\}$ if $t>3,$ and $s:=0$ otherwise; $\varepsilon_0:=\Psi(\Omega^{s+1}\theta(i_2)\coprod\Omega^{s+2}\theta(i_1)),$ $\theta_{1,2}:=\theta(i_1)\coprod\theta(i_2),$ $\theta_{0,1}:=\theta(i_0)\coprod\theta(i_1)$ and the Ext-projective cover $Q:=Q_0(\coprod_{i=1}^t\theta(i))$ in $\F(\theta).$ Then,  
\begin{itemize}
\item[(a)] if $\F(\theta)$ satisfies the 3-finitistic property then  $$\pfd\,\F(\theta)\leq s+4+\varepsilon_0+\Psi(\Omega^{s+\varepsilon_0+3}(\theta_{1,2}\coprod Q)\coprod\Omega^{s+\varepsilon_0+4}\theta_{0,1});$$ 
\item[(b)] if $\F(\theta)$ satisfies the 3-cardinal property then  $$\pfd\,\F(\theta)\leq s+2+\mathrm{max}\,(\Psi(\Omega^{s+1}\theta(i_2)\coprod\Omega^{s+2}Q),\,\Psi(\Omega^{s+1}(\theta_{1,2}\coprod Q)\coprod\Omega^{s+2}\theta_{0,1})).$$ 
\end{itemize}
\end{teo}
\begin{dem} Let $M\in\F(\theta)$ be such that  $\pd\,M$ is finite. If $\mathrm{card}\,(\infty_\theta\cap\supp\,(M))\leq 1,$ then we get, from \ref{Teofd1}, that $\pd\,M\leq s;$ proving the result in this case.  
\

Suppose that $\mathrm{card}\,(\infty_\theta\cap\supp\,(M))=2.$ Then,
we can apply \ref{Teofd2} to the stratifying system $(\theta_M,\leq);$ and so, from \ref{ProIT} (c), we get that $$\pd\,M\leq s+2+\Psi(\Omega^{s+1}\,\theta_{1,2}\coprod\Omega^{s+2}\,\theta_{0,1});$$ thus, the result is also true in this case. As can be seen, up to now, we have not used either of the "3-properties".
\

Assume that $\infty_\theta\subseteq\supp\,(M).$ We will proceed in a very similar way as we did in the proof of \ref{Teofd2}. Therefore, we have the following exact sequence in $\F(\theta)$
\begin{equation}\label{Teofd3-1} 0\rightarrow M_2\rightarrow M_{i_0}\rightarrow M^{(2)}\rightarrow 0 
\end{equation} such that $\pd\,M^{(2)}$ is finite, $\infty_\theta\subseteq\supp\,(M^{(2)}),$ $\mini\,(M^{(2)})=i_0,$ $\maxi\,(M^{(2)})=i_2$ and $\pd\,M\leq\maxi\,(s,\pd\,M^{(2)}).$ Now, we have all the ingredients we need to prove the theorem in this case. 
\

(a) Suppose that $\F(\theta)$ satisfies the 3-finitistic property. So, we take the Ext-projective cover of $M^{(2)}$ obtaining the following exact sequence in $\F(\theta)$
\begin{equation}\label{Teofd3-2} 0\rightarrow N\rightarrow
Q_0\,(M^{(2)})\rightarrow M^{(2)}\rightarrow 0, 
\end{equation} where $i_0=\mini(M^{(2)})<\mini\,(N)$ and $Q_0\,(M^{(2)})\in\add\,\coprod_{j\geq i_0}\,Q(j).$ In particular, we have that $i_0\not\in\supp\,(N).$ We consider the following three cases.
\

Case 1: $\infty_\theta\cap\supp\,(N)=\emptyset.$ Then, by \ref{LemaH2} (a), we get that  $\pd\,N\leq s;$ and so, by applying \ref{ProIT} and \ref{LemaH} to \eqref{Teofd3-2}, we  get $\pd\,M^{(2)}\leq s+2+\Psi\,(\Omega^{s+1}\,Q);$ proving that $$\pd\,M\leq s+2+\Psi\,(\Omega^{s+1}\,Q).$$
\

Case 2: $\mathrm{card}\,(\infty_\theta\cap\supp\,(N))=1.$ So, applying \ref{LemaH2} (b) to the stratifying system $(\theta_N,\leq),$ we get that $\Omega^{s+1}(N)\coprod P\simeq\Omega^{s+1}(\theta(j)^{m_j})\coprod P',$ where $P$ and $P'$ are projective $R$-modules, $j\in\supp\,(N)\cap\{i_1,i_2\}$ and $m_j:=[N:\theta(j)].$ Thus, by applying \ref{ProIT} and \ref{LemaH} to \eqref{Teofd3-2}, we  get $$\pd\,M\leq s+2+\Psi\,(\Omega^{s+1}\,(\theta_{1,2}\coprod Q)).$$
\

Case 3: $\mathrm{card}\,(\infty_\theta\cap\supp\,(N))=2.$ Hence, $i_1,i_2\in\supp\,(N);$ and so, since $\pd\,M^{(2)}<\infty,$ we get that $\pd\,N<\infty.$ Therefore, from \ref{Teofd2}, we conclude that $\pd\,N\leq s+2+\varepsilon_0.$ Thus, applying \ref{ProIT} and \ref{LemaH} to \eqref{Teofd3-2}, we  get $$\pd\,M\leq s+4+\varepsilon_0+\Psi\,(\Omega^{s+\varepsilon_0+3}\, Q).$$ 
The result follows by applying Lemma \ref{LemaH} ($\varepsilon_0+2$ times) to the first, second and third inequalities.
\

(b) Suppose that $\F(\theta)$ satisfies the 3-cardinal property. Since $\maxi\,(M^{(2)})=i_2,$ there is an exact sequence in $\F(\theta)$
\begin{equation}\label{Teofd3-3} 0\rightarrow \theta(i_2)^{m_2}\rightarrow M^{(2)}\rightarrow M^{(1)}\rightarrow 0, 
\end{equation}
where $i_0=\mini(M^{(2)})=\mini\,(M^{(1)})$ and $\infty_\theta\cap\supp\,(M^{(1)})=\{i_0,i_1\}.$ So, by taking the Ext-projective cover of $M^{(1)},$ we construct the following exact and commutative diagram in $\F(\theta)$

$$\xymatrix{ & & 0\ar[d] & 0\ar[d] & \\
&  & K\ar@{=}[r]\ar[d] & K\ar[d] & &\\
0\ar[r] & \theta(i_2)^{m_2}\ar[r]\ar@{=}[d] & E\ar[r]\ar[d] & Q_0(M^{(1)})\ar[r]\ar[d] & 0\\
0\ar[r]& \theta(i_2)^{m_2}\ar[r] & M^{(2)}\ar[r]\ar[d] & M^{(1)}\ar[d]\ar[r] & 0\\
& & 0 & 0, & }$$ where $i_0<\mini\,(K);$ in particular, $\infty_\theta\cap\supp\,(K)\subseteq\{i_1,i_2\}.$\\
Since $Q_0(M^{(1)})$ is $\theta$-projective, we have that $E\simeq \theta(i_2)^{m_2}\coprod Q_0(M^{(1)});$ and so, from \ref{ProIT}, we get that 
\begin{equation}\label{Teofd3-4} \pd\,M^{(2)}\leq 1+\Psi\,(\theta(i_2)\coprod Q\coprod K).
\end{equation} On the other hand, since $\F(\theta)$ satisfies the 3-cardinal property, we obtain that $\mathrm{card}\,(\infty_\theta\cap\supp\,(K))\leq 1.$ If $\infty_\theta\cap\supp\,(K)=\emptyset,$ then $\pd\,K\leq s;$ and hence, from \ref{sisieq} (b), we obtain that $\Omega^{s+2}\,M^{(1)}\simeq \Omega^{s+2}\,Q_0(M^{(1)}).$ Hence, by applying \ref{W} to \eqref{Teofd3-3}, we obtain, from \ref{LemaH}, the following inequality
$$\pd\,M\leq 2+s+\Psi\,(\Omega^{s+1}\,\theta(i_2)\coprod\Omega^{s+2}\,Q)$$ since $Q_0(M^{(1)})\in\add\,Q.$\\
Finally, in case $\mathrm{card}\,(\infty_\theta\cap\supp\,(K))=1,$ we conclude from \ref{LemaH2} that $\Omega^{s+1}\,(K)\coprod P\simeq \Omega^{s+1}\,(\theta(j)^{m_j})\coprod P'$ for $j\in\infty_\theta\cap\supp\,(K)$ and some projective $R$-modules $P$ and $P'.$ Then, from \eqref{Teofd3-4} and \ref{LemaH}, we get $$\pd\,M\leq 2+s+\Psi\,(\Omega^{s+1}\,(\theta_{1,2}\coprod Q)).$$ Then, joining the first and the last two inequalities, the result follows.
\end{dem}

\section{Examples}

Let $(\theta,\leq)$ be a stratifying system of size $t.$ We recall that $\infty_\theta=\{i\in[1,t]\, :\, \pd\,\theta(i)=\infty\}.$ In Section 4, we have proven, without extra conditions on
$\F(\theta),$ that $\pfd\,\F(\theta)$ is finite when $\mathrm{card}\,\infty_\theta=1,2.$  In the case $\mathrm{card}\,\infty_\theta=3,$ we were not able to prove that $\pfd\,\F(\theta)$ is finite without extra conditions. Moreover, if $\F(\theta)$ satisfies one of the ``3-properties" defined in \ref{3-properties}, we have proven in \ref{Teofd3} that $\pfd\,\F(\theta)$ is finite. Therefore, if a counter-example of the finitistic dimension conjecture exists, it might be given by a stratifying system. For example, we could start by looking for a stratifying system $(\theta,\leq)$ with $\mathrm{card}\,\infty_\theta=3$ and such that both of the ``3-properties" do not hold for $\F(\theta).$ 
\

We point out that the problem, in the proof of The Homological Stratifying System Conjecture I, appeared when $\mathrm{card}\,\infty_\theta=3.$ Note that the same kind of problems appear in at least two known results when attempting to prove, using Igusa-Todorov's function, the finitistic dimension conjecture.

\begin{ex} K. Igusa and G. Todorov in \cite{IT}.
\

The result we would like to prove: $$ \text{If }\mathrm{repdim}\,(R)<\infty\text{ then }\pfd\,(R)<\infty.$$

The result that can be proven by using Igusa-Todorov's function:
$$\text{If }\mathrm{repdim}\,(R)\leq 3\text{ then }\pfd\,(R)<\infty.$$
\end{ex}

\begin{ex} C.C. Xi in \cite{Chang}.
\

The result we would like to prove: If $I_j$ for $1\leq j\leq n$ is a family of ideals of $R$ such that $I_1I_2\cdots I_n=0$ and $R/I_j$ is of finite representation type for any $j,$ then $\pfd\,(R)<\infty.$
\

The result that can be proven by using Igusa-Todorov's function: Let $I_j$ with $1\leq j\leq n\geq 2$ be a family of ideals of $R$ such that $I_1I_2\cdots I_n=0$ and $R/I_j$ is of finite representation type for any $j.$ If $\pd\,I_j<\infty$ (as left $R$-module) and $\pd\,I_j=0$ (as right $R$-module) for $j\geq 3,$ then $\pfd\,(R)<\infty.$
\end{ex}

\begin{ex} Let $R$ be the quotient path $k$-algebra given by the following quiver

$$
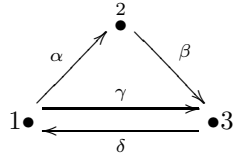
\xymatrix{ & \stackrel{2}{\bullet} \ar[dr]^{\beta} & \\
1 \bullet \ar[ru]^{\alpha} \ar@<1ex>[rr]^{\gamma}& & \ar@<1ex>[ll]^{\delta} \bullet 3}$$ and the relations $\gamma\delta\beta=\delta\gamma=\alpha\delta=\delta\beta\alpha=0.$ It is clear that $\mathrm{rad}^3R=0;$ and so, by the Corollary 7 in \cite{IT}, we know that 
\begin{equation}\label{Cap5-1} \pfd\,(R)\leq 2+\Psi\,(R/\mathrm{rad}\,(R)\coprod R/\mathrm{rad}^2(R)).
\end{equation} On the other hand, the structure of the indecomposable projective $R$-modules is $$P(1)=\begin{matrix} {} & 1 & {}\\ 2 & {} & 3\\ 3 & {}& {}\end{matrix}\qquad\quad P(2)=\begin{matrix} 2\\ 3\\ 1\end{matrix}\qquad P(3)=\begin{matrix} 3\\ 1\\ 3\end{matrix}.$$ Therefore the standard $R$-modules are ${}_R\Delta(1)=S(1),$ ${}_R\Delta(2)=S(2)$ and ${}_R\Delta(3)=P(3).$ Now, we will compute the minimal projective resolution of the simple $R$-module $S(2).$ The relevant parts of this resolution are the following exact sequences

\begin{align*}  
0\rightarrow  \begin{matrix} 3\\ 1\end{matrix}\rightarrow P(2)\rightarrow S(2)\rightarrow 0,\\ 
0\rightarrow  S(3)\rightarrow P(3)\rightarrow  \begin{matrix} 3\\ 1\end{matrix}\rightarrow 0,\\ 
0\rightarrow  \begin{matrix} 1\\ 3\end{matrix}\rightarrow P(3)\rightarrow S(3)\rightarrow 0,\\ 0\rightarrow  \begin{matrix} 2\\ 3\end{matrix}\rightarrow P(1)\rightarrow \begin{matrix} 1\\ 3\end{matrix}\rightarrow 0,\\ 
0\rightarrow  S(1)\rightarrow P(2)\rightarrow  \begin{matrix} 2\\ 3\end{matrix}\rightarrow 0,\\
0\rightarrow  S(3)\coprod\begin{matrix} 2\\ 3\end{matrix}\rightarrow P(1)\rightarrow S(1)\rightarrow 0.
\end{align*} So, we get that $\pd\,S(1)=\pd\,S(2)=\pd\,S(3)=\infty.$ Therefore, $\infty_{{}_R\Delta}=\{1, 2\},$ and we can apply \ref{Teofd2} to get that 
\begin{equation}\label{Cap5-2} \pfd\,\F({}_R\Delta)\leq 2+\mathrm{min}\,(\alpha,\beta)
\end{equation} where $\alpha=\Psi\,(\Omega\,{}_R\Delta(2)\coprod\Omega^2\,{}_R\Delta(1))$ and $\beta=\Psi(\Omega\,(Q\coprod{}_R\Delta(2))).$ We will compute those numbers, and in order to do that, we have to know the Ext-projective cover $Q=Q(1)\coprod Q(2)\coprod Q(3),$ in $\F({}_R\Delta),$ of  the $R$-module ${}_R\Delta(1)\coprod {}_R\Delta(2)\coprod {}_R\Delta(3).$ It can be proven that $Q(3)=P(3);$ and moreover, that the following exact sequences satisfy the Definition \ref{epss}
\begin{align*} 
0\rightarrow {}_R\Delta(3)\rightarrow Q(2)\rightarrow \Delta(2)\rightarrow 0,\\
0\rightarrow {}_R\Delta(2)\coprod {}_R\Delta(3)\rightarrow Q(1)\rightarrow \Delta(1)\rightarrow 0,
\end{align*} where $Q(2)=\;\begin{matrix} 3 & {} & {} & {}\\ {} & 1 & {} & 2\\ {} & {} & 3 & {}\end{matrix}\qquad$ and $\qquad Q(1)=\begin{matrix} 1 & {} & {} & {} & 3\\ {} & 2 & {} & 1 & {}\\ {} & {} & 3 & {} &  {}\end{matrix}.$\\ 
Furthermore, we have the following exact sequences 
\begin{align*} 
0\rightarrow S(3)^2\rightarrow P(1)\coprod P(3) \rightarrow Q(1)\rightarrow 0,\\
0\rightarrow \begin{matrix} 3\\ 1\end{matrix}\rightarrow P(2)\coprod P(3) \rightarrow Q(2)\rightarrow 0.
\end{align*} Then we have $M:=\Omega\,(Q\coprod{}_R\Delta(2))=S(3)^2\coprod \begin{matrix} 3\\ 1\end{matrix}\coprod\begin{matrix} 3\\ 1\end{matrix}.$ After some calculations, we can see that $\beta=\Psi\,(M)=0.$ On the other hand, $M':=\Omega\,{}_R\Delta(2)\coprod\Omega^2\,{}_R\Delta(1)=\begin{matrix} 3\\ 1\end{matrix}\coprod\begin{matrix} 1\\ 3\end{matrix}\coprod S(1);$ and so, it can be proven (using Fibonacci's sequence!) that $\alpha=\Psi\,(M')=1.$ Hence, from \eqref{Cap5-2}, we get that $\pfd\,\F({}_R\Delta)\leq 2.$ In a very similar way, it can be seen that  $\Psi\,(R/\mathrm{rad}\,(R)\coprod R/\mathrm{rad}^2(R))=2;$ hence, from \eqref{Cap5-1}, we obtain that $\pfd\,(R)\leq 4.$
\end{ex}

\begin{ex} Consider now, the following quotient path $k$-algebra $A$ given by the quiver
$$
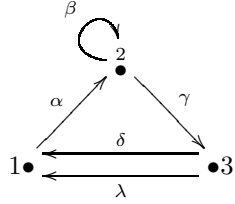
\xymatrix{ & \stackrel{2}{\bullet} \ar@(l,u)[]^{\beta} \ar[dr]^{\gamma} & \\
1 \bullet \ar[ru]^{\alpha} & & \ar@<1ex>[ll]^{\lambda} \ar@<-1ex>[ll]_{\delta}\bullet 3}$$
and the relations given in such a way that all the paths, except for $\beta^3\alpha,$ of length 4 are equal to zero and $\gamma\alpha =\gamma\beta\alpha.$ In this case, the standard $A$-modules are $${}_A\Delta(1)=S(1),\qquad {}_A\Delta(2)=\begin{matrix} 2\\ 2\\ 2\\ 2\end{matrix}\qquad\text{ and }\quad {}_A\Delta(3)=P(3).$$ It can be seen that $\pd\,{}_A\Delta(1)=\pd\,{}_A\Delta(2)=\infty;$ and so, from \ref{Teofd2}, we get that $\pfd\,\F({}_A\Delta)$ is finite. We assert that $\F({}_A\Delta)$ is of infinite representation type. Indeed, for each $t\in k,$ consider the matrices

$$M_t(\alpha):=\begin{bmatrix} 1 & 0 & 0 & 0 \\ 0 & 0 & 0 & 0 \\ 0 & 0 & 0 & 0 \\0 & 1 & 0 & 0 \\ 0 & 0 & 1 & 0 \\0 & 0 & 1 & 0 \\ 0 & 0 & 0 & 1 \\ 0 & 0 & 0 & t \end{bmatrix}\quad M_t(\beta):=\begin{bmatrix} 0 & 0 & 0 & 0 & 0 & 0 & 0 & 0 \\ 0 & 0 & 0 & 0 & 0 & 0 & 0 & 0 \\ 1 & 0 & 0 & 0 & 0 & 0 & 0 & 0 \\ 0 & 1 & 0 & 0 & 0 & 0 & 0 & 0 \\ 0 & 0 & 1 & 0 & 0 & 0 & 0 & 0 \\0 & 0 & 0 & 1 & 0 & 0 & 0 & 0 \\ 0 & 0 & 0 & 0 & 1 & 0 & 0 & 0 \\ 0 & 0 & 0 & 0 & 0 & 1 & 0 & 0\end{bmatrix}.$$ It can be proven, that the following representation $M_t$

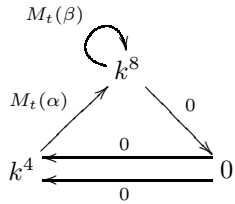
$$\xymatrix{ & k^8 \ar@(l,u)[]^{M_{t}(\beta)} \ar[dr]^{0} & \\
k^4 \ar[ru]^{M_t(\alpha)} & & \ar@<1ex>[ll]^{0} \ar@<-1ex>[ll]_{0} 0}$$ is indecomposable for any $t\in k.$ Moreover, $M_{t_1}\not\simeq M_{t_2}$ if $t_1\neq t_2;$ also $M_t\in\F({}_A\Delta)$ since there is an exact sequence $0\rightarrow {}_A\Delta(2)^2 \rightarrow M_t\rightarrow {}_A\Delta(1)^4\rightarrow 0.$
\end{ex}

\footnotesize
\vskip3mm \noindent Fran\c cois Huard:\\ Department of Mathematics,
Bishop's University.\\ Sherbrooke, Qu\'ebec, CANADA.\\ {\tt fhuard@ubishops.ca}

\vskip3mm \noindent Marcelo Lanzilotta M:\\ Centro de Matem\' atica (CMAT), \\ 
Instituto de Matem\'atica y Estad\'\i stica Rafael Laguardia (IMERL), Universidad de la Rep\'ublica.\\ Igu\'a 4225, C.P. 11400, Montevideo, URUGUAY.\\ {\tt marclan@cmat.edu.uy, marclan@fing.edu.uy}

\vskip3mm \noindent Octavio Mendoza Hern\'andez:\\ Instituto de Matem\'aticas, Universidad Nacional Aut\'onoma de M\'exico.\\ 
Circuito Exterior, Ciudad Universitaria,
C.P. 04510, M\'exico, D.F. M\'EXICO.\\ {\tt omendoza@matem.unam.mx}


\begin{thebibliography}{20}

\bibitem{ADL} I. \'Agoston, V. Dlab, E. Luk\'acs. Stratified Algebras. {\it Math. Report of Academy of Science, Canada} 20  (1988), 22-28. 
\bibitem{AHLU} I. \'Agoston, D. Happel, E. Luk\'acs, L. Unger. Finitistic dimension of standardly stratified Algebra. {\it Comm. in Algebra} 28(6)(2000) 2745-2752.
\bibitem{AR} M. Auslander, I. Reiten. Applications of Contravariantly Finite Subcategories. {\it Advances in Math.} 86 (1991) 111-152.
\bibitem{AR-Libro} M. Auslander, I. Reiten, S.O. Smalo. Representation Theory of Artin Algebras. {\it Cambridge University Press} (1995).
\bibitem{Chang} C. Xi. On the finitistic dimension conjecture I: related to representation-finite algebras. {\it Journal of pure and applied algebra} 193 (2004) nº 1-3, 287-305. 
\bibitem{Xi} C. Xi. On the finitistic dimension conjecture. {\it Advances in ring theory, World Sci. Publ., Hackensock, NJ,} (2005), 282-294.
\bibitem{DR} V. Dlab, C. M. Ringel. The Module Theoretical approach to   Quasi-hereditary algebras. {\it Repr. Theory and Related Topics, London Math. Soc. LNS} 168 (1992), 200-224.
\bibitem{ES} K. Erdmann, C. S\'aenz. On Standardly Stratified Algebras. {\it Communications in Algebra} 31(7) (2003), 3429-3446.
\bibitem{IT} K. Igusa, G. Todorov. On the finitistic global dimension conjecture for artin algebras. {\it Representation of algebras and related topics,} 201-204. Field Inst. Commun., 45. Amer. Math. Soc., Providence, RI, 2005.
\bibitem{MMS1} E. Marcos, O. Mendoza, C. S\'aenz.  Stratifying systems via relative simple modules. {\it Journal of algebra} 280 (2004) 472-487.
\bibitem{MMS2} E. Marcos, O. Mendoza, C. S\'aenz.  Stratifying systems via relative projective modules. {\it Comm. in Algebra} 33 (2005), 1559-1573.
\bibitem{MMS3} E. Marcos, O. Mendoza, C. S\'aenz.  Applications of stratifying systems to the finitistic dimension. {\it Journal of pure and applied algebra} 205 (2006), 393-411.
\bibitem{MS} O. Mendoza, C. S\'aenz.  Tilting categories with applications to stratifying systems. {\it Journal of algebra} 302 (2006), 419-449.
\bibitem{W} Y. Wang.  A note on the finitistic dimension conjecture. {\it Comm. in algebra} 22(7),419-449 (1994).
\bibitem{ZH} B. Zimmerman-Huisgen. The finitistic dimension  conjecture- a tale of 3.5 decades. in: Abelian groups and modules (Padova, 1994) 501-517. Math. Appl. 343, Kluwer Acad. Publ. Dordrecht, 1995.
\end{thebibliography}
\end{document}